\documentclass[a4paper, 11pt]{amsart}

\usepackage[mathscr]{eucal}
\usepackage{amsmath}
\usepackage{amssymb}
\usepackage{amsfonts}
\usepackage{amsthm}

\usepackage{cases}

\usepackage[mathscr]{eucal}
\usepackage{amsmath}
\usepackage{amssymb}
\usepackage{amsfonts}
\usepackage{amsthm}
\usepackage[dvips]{graphics, color}
\theoremstyle{plain}
\newtheorem{theorem}{Theorem}[section]
\newtheorem{lemma}{Lemma}[section]
\newtheorem{remark}{Remark}[section]

\numberwithin{equation}{section}

\def\<{\left<} \def\>{\right>}

\def\bea{\begin{eqnarray} }
\def\eea{\end{eqnarray} }
\def\be{\begin{equation} }
\def\ee{\end{equation} }
\def\qed{\ifhmode\unskip\nobreak\fi\ifmmode\ifinner\else\hskip5pt
\fi\fi\hbox{\hskip5 pt \vrule width4 pt height6 pt depth1.5 pt \hskip1pt }}

\begin{document}
\title[]{Real hypersurfaces in the complex projective plane attaining equality in a basic inequality}
\author[]{Toru Sasahara}
\address{Division of Mathematics, 
Hachinohe Institute of Technology, 
Hachinohe, Aomori 031-8501, Japan}
\email{sasahara@hi-tech.ac.jp}


\begin{abstract} 
We determine
non-Hopf hypersurfaces with constant mean curvature in the complex projective plane which attain equality in a basic inequality between the maximum Ricci
curvature  and the squared mean curvature.
\end{abstract}



\maketitle

 \section{Statement of the main theorem}

Let $M$ be a  real hypersurface in  the complex projective space  $\mathbb{C}P^n(4)$ of complex dimension $n$ and constant
holomorphic  sectional curvature $4$.  We denote by $J$ the almost 
complex structure of $\mathbb{C}P^n(4)$. 
The characteristic vector field on $M$ is defined by $\xi=-JN$ for a  unit normal vector field $N$.
If $\xi$ is a principal curvature vector at $p\in M$, then $M$ is said to be 
{\it Hopf} at $p$.  If $M$ is Hopf at every point, then $M$ is called a {\it Hopf hypersurface}.
Let $\mathcal{H}$ be the holomorphic distribution defined by $\mathcal{H}=\bigcup_{p\in M}\{X\in T_pM\ |\ \<X, \xi\>=0\}$.
If $\mathcal{H}$ is integrable and each leaf of its maximal integral manifolds is 
locally congruent to
$\mathbb{C}P^{n-1}(4)$, then $M$ is called a {\it ruled real hypersurface}, which is a
 typical example of a non-Hopf hypersurface.

For a Riemannian manifold $M$, 
let  $\overline{Ric}$ denote the maximum Ricci curvature function on $M$ defined by
\be\overline{Ric}(p)=\max\{S(X, X)\ |\ X\in T_pM, \ ||X||=1\},\nonumber\ee
where $S$ is the Ricci tensor. 
On the other hand, 
the $\delta$-invariant $\delta(2)$ of $M$ is 
defined by $\delta(2)(p)=(1/2)\tau(p)-\min\{K(\pi)\ |\ \pi\  \text{is a plane in}\ T_pM\}$,
where $\tau$ is the scalar curvature of $M$ and $K(\pi)$ is the sectional curvature of $\pi$. 
 (For general $\delta$-invariants,
see \cite{chen6} for details.)
In the case of $\dim M=3$, we have
 $\overline{Ric}(p)=\delta(2)(p)$.
%
Thus, according to Corollary 7 and Theorem 8 in \cite{chen}, 
for real hypersurfaces in  $\mathbb{C}P^2(4)$ we have the following:
\begin{theorem}[\cite{chen}]
Let  $M$ be a real hypersurface in 
${\mathbb C}P^2(4)$.
 Then we have
 \be
 \overline{Ric}\leq \frac{9}{4}\|H\|^2+5.\label{ideal}
 \ee
 If $M$ is a Hopf  hypersurface, then the equality in $(\ref{ideal})$ holds identically
 if and only if one of the following two cases occurs{\rm :}
 
 {\rm (1)} $M$ is an open portion of a geodesic sphere with radius $\pi/4${\rm ;}

{\rm (2)} $M$ is an open portion of a tubular hypersurface over a
complex quadric curve $Q_1$
 with radius $r=\tan^{-1}((1+\sqrt{5}-\sqrt{2+2\sqrt{5}})/2)=0.33311971\cdots$. 
 \end{theorem}

The next step is to classify non-Hopf hypersurfaces in ${\mathbb C}P^2(4)$ which satisfy
the equality in $(\ref{ideal})$ identically.
The main theorem of this paper is the following.
 
\begin{theorem}\label{main}
Let $M$ be a real hypersurface in $\mathbb{C}P^2(4)$ which is non-Hopf at 
every point. 
Assume that $M$ has constant mean curvature.
 Then $M$ satisfies the equality case of $(\ref{ideal})$ identically
  if and only if it is a minimal ruled real hypersurface which is  given by
  $\varpi\circ z$, where  $\varpi: S^{5}\rightarrow \mathbb{C}P^2(4)$ is the Hopf fibration and 
\be
z(u, v, \theta, \psi)=e^{\sqrt{-1}\psi}\bigl(\cos{u}\cos{v}, \cos{u}\sin{v}, (\sin{u})e^{\sqrt{-1}\theta}\bigr)\nonumber\ee
for
$-\pi/2<u<\pi/2$, \ \  $0\leq v,  \theta, \psi<2\pi$.
\end{theorem}


 \begin{remark}
 {\rm Let $M$ be a real hypersurface in the complex hyperbolic plane 
 $\mathbb{C}H^2(-4)$ of constant
 holomorphic sectional curvature $-4$. Then we have
 \be
  \overline{Ric}\leq \frac{9}{4}\|H\|^2-2. \nonumber
 \ee
 The equality sign of the inequality holds identically if and only if
 $M$ is an open part of the horosphere (see \cite{chen} or \cite{chen5}). }
 \end{remark}

\section{Preliminaries}

Let $M$ be a   real hypersurface in  $\mathbb{C}P^n(4)$.
Let us denote by $\nabla$ and
 $\tilde\nabla$ the Levi-Civita connections on $M$ and $\mathbb{C}P^n(4)$, respectively. The
 Gauss and Weingarten formulas are respectively given by
\be
 \begin{split}
 \tilde \nabla_XY&= \nabla_XY+\<AX, Y\>N, \label{gawe}\\
 \tilde\nabla_X N&= -AX
 \end{split}\nonumber
\ee
 for tangent vector fields $X$, $Y$ and a unit normal vector field $N$,
 where $A$ is the shape operator.
The mean curvature vector field $H$ is defined by 
$H=({\rm Tr}A/(2n-1))N.$
 The function ${\rm Tr}A/(2n-1)$ is called the  {\it mean curvature}.
 If it vanishes identically, then $M$ is called a {\it minimal hypersurface}.
 
For any  vector field  $X$ tangent to $M$,  we denote the tangential component of $JX$ by $PX$.
Then by the Gauss and  Weingarten formulas, we have
\be
\nabla_{X}\xi=PAX. \label{PA}
\ee

 
 We denote by $R$ the Riemannian curvature tensor of $M$. Then,
 the equations of Gauss  and Codazzi are respectively given by
 \begin{align}
 &R(X, Y)Z=\<Y, Z\>X-\<X, Z\>Y+\<PY, Z\>PX
 -\<PX, Z\>PY  \label{ga}\\
 &\hskip60pt -2\<PX, Y\>PZ
   +\<AY, Z\>AX-\<AX, Z\>AY,\nonumber\\
  &(\nabla_XA)Y-(\nabla_YA)X=\<X, \xi\>PY-\<Y, \xi\>PX-2\<PX, Y\>\xi.\label{co}
 \end{align}

 We need  the following two lemmas for later use.
 \begin{lemma}[\cite{chen}]\label{lem1}
  Let $M$ be a real  hypersurface 
   in $\mathbb{C}P^2(4)$. Then the equality sign in $(\ref{ideal})$ holds at a point $p\in M$ 
if and only if there exists an orthonormal basis $\{e_1, e_2, e_3\}$ at $p$ such that

${\rm (1)}$ $\<Pe_1, e_2\>=0$,  

${\rm (2)}$ the shape operator of $M$ in ${\mathbb C}P^2(4)$ at $p$ satisfies
\be
A= \left(
    \begin{array}{ccc}
      \alpha & \beta & 0 \\
      \beta & \gamma & 0\\ 
      0 & 0 & \mu \\
   \end{array}
  \right),\label{A} \ \ \alpha+\gamma=\mu.
\ee
\end{lemma}

\begin{lemma}[\cite{ni}]\label{lem2}
Let $M$ a  real hypersurface $M$ in $\mathbb{C}P^n(4)$ with $n\geq 2$.
We define differentiable functions $\alpha$, $\beta$ on $M$ by $\alpha=\<A\xi, \xi\>$
and $\beta=\|A\xi-\alpha\xi\|$.
Then $M$ is ruled if and only if 
the following two conditions hold{\rm :} 




{\rm (1)} the set $M_1=\{p\in M \ | \ \beta(p)\ne 0\}$ is an open dense subset of $M${\rm ;}

{\rm (2)} there is a unit vector field $U$ on $M_1$, which is orthogonal to $\xi$ and satisfies
\be 
A\xi=\alpha\xi+\beta U, \ \ AU=\beta\xi, \ \ AX=0 \label{ruled}
\ee
 for an arbitrary tangent  vector $X$ orthogonal to both $\xi$ and $U$.

\end{lemma}

\section {Proof of the main theorem}

 Let $M$ be a real hypersurface  in $\mathbb{C}P^2(4)$ which 
 is non-Hopf at every point.
Assume that $M$ has constant mean curvature 
and satisfies
   the equality case of $(\ref{ideal})$ identically.

   Let $\{e_1, e_2, e_3\}$ be a local orthonormal frame field described in Lemma \ref{lem1}. 
It follows from (1) of Lemma \ref{lem1}  that $\xi$ lies in Span$\{e_1, e_2\}$. 
Thus, we may assume that $e_1=\xi$ and $Pe_2=e_3$. Then, 
(\ref{A}) can be rewritten as
\begin{align} 
& A\xi=(\mu-\gamma)\xi+\beta e_2,\ \  Ae_2=\gamma e_2+\beta\xi, \ \ Ae_3=\mu e_3. \label{A2}
 \end{align}
Since $\xi$ is not a principal vector everywhere, 
we have $\beta\ne 0$ on $M$. The constancy of the mean curvature implies that $\mu$ is constant.
By (\ref{PA}) and (\ref{A2}), we have 
\be
\nabla_{e_2}\xi=\gamma e_3, \ \ \nabla_{e_3}\xi=-\mu e_2, \ \ \nabla_{\xi}\xi=\beta e_3. \label{n1}
\ee
Since $\<\nabla e_i, e_j\>=-\<\nabla e_j, e_i\>$ holds, it follows from (\ref{n1}) that
\begin{align}
& \nabla_{e_2}e_2 =\kappa_1e_3, \ \ \nabla_{e_3}e_2=\kappa_2e_3+\mu\xi, \ \ 
 \nabla_{\xi}e_2 =\kappa_3e_3, \nonumber \\
& \nabla_{e_2}e_3 =-\kappa_1e_2-\gamma\xi, \ \ \nabla_{e_3}e_3=-\kappa_2e_2, \ \ \nabla_{\xi}e_3=
-\kappa_3e_2-\beta\xi\nonumber
\end{align}
for some smooth functions  $\kappa_1$, $\kappa_2$ and $\kappa_3 $.
 %

With respect to the Gauss-Codazzi equations, we  are going to 
state only equations that are useful in this  proof.

From the equation  (\ref{co}) of Codazzi,  we obtain:
 \begin{itemize}
\item  for $X=e_2$ and  $Y=\xi$, by comparing the coefficient of $e_3$,
\be
 \beta\kappa_1+(\mu-\gamma)\kappa_3=\beta^2+\gamma^2-1.\label{eq1}
\ee

\item for $X=e_3$ and $Y=\xi$, by noting that $\mu$ is constant and $\beta\ne 0$, 
\begin{align}
e_3\beta &=\mu^2-2\mu\gamma-\kappa_3(\mu-\gamma)+\beta^2+1, \label{eq3}\\
e_3\gamma&=2\beta\mu+\beta\gamma-\beta k_3,\label{cd5}\\
 \kappa_2 &=0.\label{cd4}
\end{align}
\item for $X=e_2$ and  $Y=e_3$, by comparing the coefficient of $e_2$, 
\be
e_3\gamma=-\mu\kappa_1+\kappa_1\gamma+\beta\gamma+2\beta\mu.\label{cd7}
\ee
 \end{itemize}
Eliminating $e_3\gamma$ from (\ref{cd5}) and (\ref{cd7}) yields
 \be
 (\mu-\gamma)\kappa_1-\beta\kappa_3=0.
 \label{eq4}
 \ee 
  By solving (\ref{eq1}) and (\ref{eq4}) for $\kappa_1$ and $\kappa_3$, we get
 \begin{align}
 \kappa_1 &=\frac{\beta(\beta^2+\gamma^2-1)}{(\mu-\gamma)^2+\beta^2}, \label{k1}\\
 \kappa_3 &=\frac{(\mu-\gamma)(\beta^2+\gamma^2-1)}{(\mu-\gamma)^2+\beta^2}.\label{k3}
 \end{align}
 
 By applying the equation (\ref{ga}) of Gauss for $X=e_2$, $Y=Z=e_3$,  comparing the 
 coefficient of $e_2$ and using (\ref{cd4}),
 we obtain 
 \be
 e_3\kappa_1-2\mu\gamma-\kappa_1^2-(\gamma+\mu)\kappa_3-4=0.\label{gauss}
 \ee
 Substituting (\ref{k1}) and (\ref{k3}) into (\ref{eq3}), (\ref{cd7}) and (\ref{gauss}), we obtain
 \begin{align}
 &e_3\beta=(\mu-2\gamma)\mu+\beta^2+1-\frac{(\mu-\gamma)^2(\beta^2+\gamma^2-1)}{(\mu-\gamma)^2+\beta^2}, \label{b}\\
& e_3\gamma=\beta(\gamma+2\mu)+\frac{(\gamma-\mu)\beta(\beta^2+\gamma^2-1)}{(\mu-\gamma)^2+\beta^2},\label{r}\\
&\{(3\beta^2+\gamma^2-1)((\mu-\gamma)^2+\beta^2)-2\beta^2(\beta^2+\gamma^2-1)\}e_3\beta  \label{g1}\\
 &+\{2\beta\gamma((\mu-\gamma)^2+\beta^2)+2(\mu-\gamma)\beta(\beta^2+\gamma^2-1)\}e_3\gamma\nonumber \\
 &-2\mu\gamma\{(\mu-\gamma)^2+\beta^2\}^2-\beta^2(\beta^2+\gamma^2-1)^2\nonumber\\
 &+(\gamma^2-\mu^2)(\beta^2+\gamma^2-1)\{(\mu-\gamma)^2+\beta^2\}-4\{(\mu-\gamma)^2+\beta^2\}^2=0.\nonumber
 \end{align}
 Substituting   (\ref{b}) and   (\ref{r}) into   (\ref{g1}) gives 
 \be
 (\mu-\gamma)f(\beta, \gamma)=0, \label{f}\ee
where  $f(\beta, \gamma)$ is given by the following polynomial.
\be
\begin{split}
f(\beta, \gamma):=& 2\mu\gamma^4-(4\mu^2-1)\gamma^3+(3\mu^2+4\beta^2-6)\mu\gamma^2-\{\mu^4+(4\beta^2-7)\mu^2-\beta^2-1\}\gamma\nonumber\\
&+(\beta^2-2)\mu^3+(2\beta^4-2\beta^2-1)\mu.\nonumber
\end{split}
\ee
 
 By (\ref{f}), our discussion is divided into two cases.
 
 \medskip
 {\it Case} (a):\ \ $\mu-\gamma=0$.\ \  In this case,  from (\ref{eq4}) we obtain $\kappa_3=0$. 
 Therefore, by (\ref{cd5})  and the constancy of $\mu$, we see
  that $\mu=\gamma=0$.

 \medskip
{\it Case} (b):\ \ $f(\beta, \gamma)=0$.\ \ In this case, differentiating $f(\beta, \gamma)=0$ along $e_3$, by using (\ref{b}) and (\ref{r}) we obtain
\be
\begin{split}
& 8\mu\gamma^6-(24\mu^2-4)\gamma^5+(30\mu^2+24\beta^2-15)\mu\gamma^4\\
 &-\{20\mu^4+(48\beta^2+3)\mu^2-8b^2-3\}\gamma^3\\
 &+\{7\mu^5+(36\beta^2+45)\mu^3+(24\beta^4-10\beta^2-2)\mu\}\gamma^2\\
&-\{\mu^6+(12\beta^2+44)\mu^4+(24\beta^4+19\beta^2+2)\mu^2-4\beta^4-3\beta^2+1\}\gamma\\
 &+(\beta^2+13)\mu^5+(6\beta^4+19\beta^2+1)\mu^3+(8\beta^6+5\beta^4-2\beta^2+1)\mu=0. \label{f2}
\end{split}
\ee
 The resultant of $f(\beta, \gamma)$ and the LHS of (\ref{f2}) with respect to $\gamma$ is given by
  \be
  202500(\mu^2-1)^4\beta^4\mu^6\{4\mu^2\beta^2+(\mu^2-1)^2\}^2.\label{re}
  \ee
Since $\beta\ne 0$, we have $4\mu^2\beta^2+(\mu^2-1)^2\ne 0$.
By changing the sign of $N$ if necessary,
we may assume that $\mu\geq 0$. Thus, from (\ref{re}) we get
$\mu\in\{0, 1\}$. 
If $\mu=1$, then  equations $f(\beta, \gamma)=0$ and (\ref{f2}) can be reduced to 
\begin{align*}
& 2\beta^2+2\gamma^2+\gamma-3=0,\\
& 8\beta^4+(16\gamma^2-4\gamma+3)\beta^2+(\gamma-1)^2(8\gamma^2+12\gamma+15)=0,
\end{align*}
 respectively. This system of  equations has a unique solution $(\beta, \gamma)=(0, 1)$,
 which contradicts $\beta\ne 0$.
 Hence, we have $\mu=0$. Then, equation $f(\beta, \gamma)=0$ becomes 
 $\gamma(\beta^2+\gamma^2+1)=0$, which shows that
 $\gamma=\alpha=0$. 
 
 \medskip
 From the above argument, the shape operator satisfies
 \be
 A\xi=\beta e_2, \ \ 
Ae_2=\beta\xi, \ \  Ae_3=0\nonumber 
 \ee 
at each point,  where $\beta\ne 0$. 
 By Lemma \ref{lem2} we conclude that 
 $M$ is a minimal ruled real hypersurface.
 According to \cite{ada}, it is congruent to the real hypersurface described in Theorem \ref{main}. 
 
 The converse is clear from Lemma \ref{lem1} and (\ref{ruled}).
The proof is finished.
\qed

 \end{document}